\documentclass[12pt]{article}

\usepackage{amsmath,amssymb,amsfonts, amsthm}

\newcommand{\Ann}{\mbox{Ann}\,}
\newcommand{\Hom}{\mbox{Hom}\,}

\newcommand{\Spec}{\mbox{Spec}\,}

\newcommand{\Supp}{\mbox{Supp}\,}

\newcommand{\depth}{\mbox{depth}\,}
\renewcommand{\dim}{\mbox{dim}\,}
\newcommand{\cd}{\mbox{cd}\,}
\newcommand{\amp}{\mbox{amp}\,}
\newcommand{\Min}{\mbox{Min}}
\newcommand{\pd}{\mbox{proj.dim}\,}

\newcommand{\h}{\mbox{ht}\,}
\newcommand{\E}{\mbox{E}}
\renewcommand{\H}{\mbox{H}}
\newcommand{\V}{\mbox{V}}

\newcommand{\G}{\Gamma}

\newcommand{\T}{\mathrm}

\newcommand{\lo}{\longrightarrow}

\newcommand{\fa}{\frak{a}}
\newcommand{\fb}{\frak{b}}

\newcommand{\fp}{\frak{p}}
\newcommand{\fq}{\frak{q}}

\date{}
\begin{document}

\title{Cohomological dimension of complexes
}

\author{Mohammad T. Dibaei $^{\it ab}$ and Siamak Yassemi$^{\it
ca}$\\
{\small\it $(a)$ Institute for Studies in Theoretical Physics and
Mathematics}\\
{\small\it $(b)$ Department of Mathematics, Teacher Training
University}\\
{\small\it $(c)$ Department of Mathematics, University of Tehran.}
}
\maketitle

\begin{abstract}

\noindent In the derived category of the category of modules over
a commutative Noetherian ring $R$, we define, for an ideal $\fa$
of $R$, two different types of cohomological dimensions of a
complex $X$ in a certain subcategory of the derived category,
namely $\cd(\fa, X)=\sup\{\cd(\fa, \H_{\ell}(X))-\ell|\ell\in\Bbb
Z\}$ and $-\inf{\mathbf R}\G_{\fa}(X)$, where $\cd(\fa,
M)=\sup\{\ell\in\Bbb Z|\H^{\ell}_{\fa}(M)\neq 0\}$ for an
$R$--module $M$. In this paper, it is shown, among other things,
that, for any complex $X$ bounded to the left, $-\inf {\mathbf
R}\G_{\fa}(X)\le\cd(\fa, X)$ and equality holds if indeed $\H(X)$
is finitely generated.

\end{abstract}

\vspace{.2in}

\baselineskip=18pt

\vspace{.3in}

\noindent{\bf 0. Introduction}

\vspace{.2in}

Let $R$ be a commutative Noetherian ring of finite dimension $d$
and $\fa$ be an ideal of $R$. For an $R$--module $M$,
$\G_{\fa}(M)$ is defined to be the submodule of $M$ consisting of
all elements of $M$ which are vanished by some power of $\fa$. It
has been an interesting question to know when $\H^i_{\fa}(M)$, the
i--th right derived functor of $\G_{\fa}(-)$ applied on $M$, is
zero (see [{\bf Hu}]). The cohomological dimension of $M$ with
respect to $\fa$ is defined as
$$\cd(\fa, M)=\sup\{i\in\Bbb Z|\H^i_{\fa}(M)\neq 0\}.$$
In [{\bf G}] Grothendieck has shown that $\cd(\fa, M)$ has a lower
bound and an upper bound $\depth M$ and $\dim M$ respectively. The
cohomological dimension has been studied by several authors. In
[{\bf Fa}] Falting and in [{\bf HL}] Huneke--Lyubeznik have found
several upper bounds for cohomological dimension. In [{\bf DNT}],
some lower bounds have been obtained.

In section 1, we study the cohomological dimension of a module in
more details. We show that for an $R$--module $M$, $\cd(\fa,
M)\le\cd(\fa, R/\fp)$ for some $\fp\in\Supp_R( M)$ and equality
holds if $M$ is finite (that is finitely generated) and in this
case $\fp$ may be taken from the set of minimal elements of
$\Supp_R(M)$ (Theorem 1.3).

In the derived category ${\cal D}(R)$, for a complex $X$, the
$\fa$--depth of $X$ is defined by the following formula, cf. [{\bf
I}; Section 2]:
$$\depth(\fa,X)=\inf\{\ell\in\Bbb Z|\H_{-\ell}({\mathbf
R}\Hom_R(R/\fa,X))\neq 0\}.$$ In [{\bf FI}; Theorem 2.1], Foxby
and Iyengar have shown that, for any complex $X$, $$\depth(\fa,
X)=-\sup{\mathbf R}\Gamma_\fa(X).$$

In this paper we are specially interested in $\inf{\mathbf
R}\Gamma_\fa(X)$. It is clear that for an $R$--module $M$,
$-\inf{\mathbf R }\Gamma_\fa(M)=\cd(\fa, M)$ and so by
Grothendieck's result, we have $-\inf{\mathbf
R}\Gamma_\fa(M)\le\dim(M)$. What is cited is our motivation to
study the invariant $-\inf{\mathbf R}\Gamma_\fa(X)$. On the other
hand, for a complex $X$ bounded to the left, in consistent with
the definition of the dimension of $X$ as $\dim X =
\sup\{\dim_R\H_{\ell}(X)-\ell|\ell\in\Bbb Z \}$, we may define
cohomological dimension of $X$ as
$\cd(\fa,X)=\sup\{\cd(\fa,\H_{\ell}(X))-\ell|\ell\in\Bbb Z\}$. Our
purpose is to clarify the relationship between $\inf{\mathbf
R}\Gamma_\fa(X)$ and $\cd(\fa, X)$.

We first study $\cd(\fa,X)$ in section 2. As $\cd(\fa, M)$ is well
behaved when $M$ is a finite module, most results of this section
are in the case $X$ has finite homology modules, for example: If
$X,~Y$ are complexes bounded to the left with finite homology
modules then
$$\cd(\fa, X\otimes^{\mathbf L}_RY)=\sup\{\cd(\fa,
\H_{t}(X)\otimes_R\H_{\ell}(Y))-\ell-t|t, \ell\in\Bbb
Z\}\,\,\,\,\mbox{(see Theorem 2.9),}$$

\noindent and we have
$$\cd(\fa, X\otimes^{\mathbf L}_RY)\le\cd(\fa, X)-\inf Y\,\,\,\,\mbox{(see Theorem 2.10)}.$$

Moreover, there are some results which $X$ does not need to have
finite homology modules, for example: If $X$ is a bounded to the
left complex, then
$$\cd(\fa, X)\le\cd(\fa, R)-\inf X\,\,\,\,\mbox{(see Proposition
2.5)}.$$

In section 3, we compare $\cd(\fa, X)$ with $-\inf{\mathbf
R}\Gamma_\fa(X)$. More precisely, we show that for any bounded to
the left complex $X$ the invariant $-\inf{\mathbf R}\Gamma_\fa(X)$
has a lower and an upper bounds $\depth(\fa, X)$ and $\cd(\fa, X)$
respectively; and it takes its upper bound when indeed the
homology modules of $X$ are finite (see Theorem 3.2 and Theorem
3.3).

\vspace{.3in}

\noindent{\bf 1. Cohomological dimension of a module}

\vspace{.2in}

In this section, we have a brief look at the cohomological
dimension of a module to find some extra properties (see Theorem
1.2) and to give an extension of [{\bf DNT}; Theorem 2.2] (see
Theorem 1.5).

First recall the well--known fact about $\cd(\fa, R)$ which states
that

(1.0.1)\,\,\,\,\, For any $R$--module $M$, $\cd(\fa, M)\le\cd(\fa,
R)$.

The following result will be helpful to remove the finiteness
condition on the modules.

\vspace{.2in}

\noindent{\bf Theorem 1.1.} If $M$ is an $R$--module, then
$$\cd(\fa, M)\le\sup\{\cd(\fa,N)| N\,\,\,\T{is\,\, a\,\, finite\,\,
submodule\,\, of}\,\,\, M\}.$$

\vspace{.1in}

\noindent{\it Proof.} It is well--known that $M$ is equal to the
direct limit of its finite submodules. Now the assertion follows
from the fact that the local cohomology functor commutes with the
direct limit.\hfill$\square$

\vspace{.2in}

\noindent{\bf Lemma 1.2.} If $M$ is an $R$--module and $\fb$ is an
ideal with $\fb\subseteq\Ann_R(M)$, then $$\cd(\fa, M)\le\cd(\fa,
R/\fb).$$

\vspace{.1in}

\noindent{\it Proof.} As $M$ is an $R/\fb$--module, we have
\[ \begin{array}{rl} \cd(\fa,
M) &=\cd(\fa (R/\fb), M)\\
 & \le\cd(\fa (R/\fb), R/\fb)\\
 & =\cd(\fa, R/\fb).
\end{array} \]
\noindent The equalities hold by independence theorem on local
cohomology and the inequality holds by (1.0.1).\hfill$\square$

\vspace{.2in}

\noindent{\bf Theorem 1.3.} If $M$ is an $R$--module with finite
cohomological dimension with respect to $\fa$, then $\cd(\fa,
M)\le\cd(\fa, R/\frak p)$ for some $\fp\in\Supp_R(M)$. Moreover,
if $M$ is finite, the equality holds and $\fp$ can be taken from
the set, $\Min\Supp_R(M)$, of minimal elements of $\Supp_R(M)$.

\vspace{.1in}

\noindent{\it Proof.} By Theorem 1.1, there exists a finite
submodule $K$ of $M$ such that $\cd(\fa, M)\le\cd(\fa, K)$, so it
is enough to show the assertion for $K$. Assume the contrary.
There is a chain $0=K_0\subset K_1\subset\cdots\subset K_n=K$ of
submodules of $K$ such that, for each $i$, $K_i/K_{i-1}\cong
R/\fp_i$, where $\fp_i\in\Supp_R(K)$. Set $t=\cd(\fa, K)$, so we
have $\H^t_{\fa}(R/\fp_i)=0$ for $1\le i\le n$. Thus from the
exact sequences $\H^t_\fa(K_{i-1})\to\H^t_\fa(K_i)\to 0$,
$i=1,2,\ldots, n$, we eventually get $\cd(\fa, K_1)\ge t$ which is
a contradiction.

Let $M$ be finite and so $\cd(\fa, M)\le\cd(\fa, R/\fp)$ for some
$\fp\in\Supp_R(M)$. On the other hand $\cd(\fa, R/\fp)\le\cd(\fa,
M)$, c.f. [{\bf DNT}; Theorem 2.2]. Hence the equality holds.

If $\fp$ does not belong to the set $\Min\Supp_R(M)$ then there
exists $\fq\in\Min\Supp_R(M)$ with $\fq\subset\fp$. Now by using
Lemma 1.2, we have that $\cd(\fa, R/\fp)\le\cd(\fa, R/\fq)$, that
is $\cd(\fa, R/\fq)=\cd(\fa, M)$, and the claim
follows.\hfill$\square$

\vspace{.2in}

The following result is a generalization of [{\bf DNT}; Theorem
2.2].

\vspace{.1in}

\noindent{\bf Theorem 1.4.} Let $N$ and $M$ be $R$--modules and
$M$ finite. If $\Supp_R(N)\subseteq\Supp_R(M)$, then $\cd(\fa,
N)\le\cd(\fa, M)$.

\vspace{.1in}

\noindent{\it Proof.} By Theorem 1.1, $\cd(\fa, N)\le\cd(\fa, K)$
for some finite submodule $K$ of $N$. Since
$\Supp_R(K)\subseteq\Supp_R(M)$, we have by [{\bf DNT}; Theorem
2.2], that $\cd(\fa, K)\le\cd(\fa, M)$.\hfill$\square$

\vspace{.2in}

\noindent{\bf Corollary 1.5.} Let $\varphi\,:R\to S$ be a ring
homomorphism and let $M$ be a finite $R$--module. Then
$$\cd_S(\fa S,S\otimes_RM)\le\cd_R(\fa,M).$$
If $S$ is faithfully flat then the equality holds.

\vspace{.1in}

\noindent{\it Proof.} By independence theorem
$\H^i_{\fa}(S\otimes_R M)\cong\H^i_{\fa S}(S\otimes_RM)$ as
$R$--modules. As $\Supp_R(S\otimes_RM)\subseteq\Supp_R(M)$, we
have, by Theorem 1.1, $\cd_S(\fa S, S\otimes_RM)\le\cd_R(\fa, M)$.
The final claim is clear.\hfill$\square$

\vspace{.2in}

The following example shows that the finiteness condition on $M$
is not redundant in Theorem 1.4 and second part of Theorem 1.3.

\vspace{.1in}

\noindent{\bf Example 1.6.} Choose a ring $R$, a prime ideal $\fp$
and an ideal $\fa$ such that $\fa\nsubseteq\fp$. We have
$\Supp_R(R/\fp)\subseteq\Supp_R(\E(R/\fp))$, where $\E(R/\fp)$ is
the injective envelope of $R/\fp$ as $R$--module. We observe that
$\cd(\fa, \E(R/\fp))=0$ but $0<\cd(\fa, R/\fp)$.

\vspace{.3in}

\noindent{\bf 2. Cohomological dimension of a complex}

\vspace{.2in}

An $R$-complex $X$ is a sequence of $R$-modules $X_\ell$ and
$R$-linear maps $\partial_\ell^X, \ell\in \mathbb{Z}$,
$$X=\cdots \lo X_{\ell+1} \stackrel{\partial^X_{\ell+1}}{\lo} X_\ell
\stackrel{\partial_\ell^X}{\lo} X_{\ell-1} \lo \cdots. $$ The
module $X_\ell$ is called the module in degree $\ell$, and the map
$\partial_\ell^X:X_\ell\lo X_{\ell-1}$ is the $\ell$-th
differential, and $\partial_\ell^X \partial_{\ell+1}^X=0$ for all
$\ell\in \mathbb{Z}$. An $R$--module $M$ is thought of as the
complex $M=0\lo M\lo 0$,with $M$ in degree zero.

The {\em supremum} and  {\em infimum} of $X$ are defined by
\[ \begin{array}{rl} \sup \: X & =
\sup \: \{ {\ell} \in \Bbb Z | \H_{ \ell } (X) \neq 0 \}
\\ [.1in] \inf \: X & = \inf \: \{ {\ell} \in \Bbb Z | \H_{
\ell } (X) \neq 0 \}
\end{array} \]

\noindent Denote $\sup X=-\infty$ and $\inf X=\infty$ if
$\H_{\ell}(X)=0$ for all $\ell$.

A morphism $\alpha: X\lo Y$ is said to be a quasi-isomorphism if
the induced morphism $\T{H}(\alpha):\T{H}(X)\lo \T{H}(Y)$ is an
isomorphism.

The {\it derived category} ${\cal D}(R)$ of the category of
$R$--complexes is the category of $R$--complexes localized at the
class of all quasi--isomorphisms. The full subcategories ${\cal
D}_{+} (R)$, ${\cal D}_{-} (R)$, ${\cal D}_{b} (R)$, and ${\cal
D}_0(R)$ consist of complexes $X$ with $\H_{\ell}(X)=0$ for,
respectively, $\ell\ll 0, ~\ell\gg 0, |\ell|\gg 0$, and $\ell\neq
0$. By ${\cal D}^f(R)$ we mean the full subcategory of ${\cal
D}(R)$ consisting of complexes $X$ with $\H_{\ell}(X)$ a finite
$R$--module for all $\ell$.

The left derived functor of the tensor product functor of
$R$-complexes is denoted by $-\otimes_R^{\T{\mathbf L}}-$, and
$\T{\mathbf R}\T{Hom}_R(-,-)$ denotes the right derived functor of
the homomorphism functor of complexes. We need the next two
inequalities for $X,~Y\in{\cal D}_+(R)$ and $Z\in{\cal D}(R)$.

(2.0.1)\,\,\,\,\,\,\,\, $\inf(X\otimes^{\mathbf L}_RY)\ge \inf
X+\inf Y$\,\,\,\,\, and

(2.0.2)\,\,\,\,\,\,\, $\sup({\mathbf R}\Hom_R(X,Z))\le\sup Z-\inf
X$.

\noindent For a complex $X$, the dimension of $X$ is defined by

(2.0.3)\,\,\,\,\,\,\, $\dim_RX=\sup\{\dim R/\fp-\inf
X_{\fp}|\fp\in\Spec R\}$.

\noindent It is shown, in [{\bf Fo1}; 16.9], that
$$\dim X = \sup\{\dim_R\H_{\ell}(X)-\ell|\ell\in\Bbb Z \}.$$

\noindent Therefore it is natural to give the following
definition:

\vspace{.1in}

\noindent{\bf Definition 2.1.} For a complex $X\in{\cal D}(R)$,
the $\fa$--cohomological dimension of $X$ is defined by

$$\cd(\fa,X)=\sup\{\cd(\fa,\H_{\ell}(X))-\ell|\ell\in\Bbb Z\}.$$

\noindent For an $R$--module $M$, this notion agrees with the
classical one. Note that $\cd(\fa, X)=-\infty$ if and only if $X$
is homologically trivial. If $\inf X=-\infty$ then
$\cd(\fa,X)=\infty$.

To find some extra information about $\cd(\fa,X)$, we also review
the notion of the height of an ideal $\fa$. If $M$ is an
$R$--module, the $M$--height of $\fa$, denoted by $\h(\fa, M)$, is
defined to be the supremum length of chains
$\fp_0\subset\fp_1\subset\cdots\subset\fp_n$ of elements of
$\Supp_RM$ with $\fp_n$ is minimal over $\fa$. Thus, when $M$ is
finite, we may write $\h(\fa, M)=\h(\fa, R/\fp)$ for some
$\fp\in\Supp_RM$, and that
$$\h(\fa, M)=\sup\{\h(\fa, R/\fp)|\fp\in\Supp_RM\}.$$
If $M$ is zero module then $\h(\fa, M)=-\infty$ by convention.

\vspace{.1in}

\noindent{\bf Definition 2.2.} For $X\in{\cal D}_+(R)$, the
$X$--height of $\fa$ is defined by
$$\h(\fa, X)=\sup\{\h(\fa, \H_{\ell}(X))-\ell|\ell\in\Bbb Z\}.$$

Now we have the following result.

\vspace{.1in}

\noindent{\bf Lemma 2.3.} If $X\in{\cal D}^f_+(R)$, then
$$\h(\fa, X)=\sup\{\h(\fa, R/\fp)-\inf X_{\fp}|\fp\in\Spec
R\}.$$

\vspace{.1in}

\noindent{\it Proof.} Consider $\ell\in\Bbb Z$ such that
$\H_{\ell}(X)\neq 0$. There exists $\fp\in\Supp_R\H_{\ell}(X)$
with $\h(\fa, R/\fp)=\h(\fa, \H_{\ell}(X))$. Thus $\h(\fa,
\H_{\ell}(X))-\ell\le\h(\fa, R/\fp)-\inf X_{\fp}$ which gives a
one side inequality.

For the other side, assume $\fp\in\Supp_RX$ and that $\ell=\inf
X_{\fp}$. Thus $\h(\fa, R/\fp)\le\h(\fa,\H_{\ell}(X))$. Now the
assertion holds.\hfill$\square$

\vspace{.2in}

The following proposition compares the invariants $\cd(\fa,
X),~\cd(\fa, R),~ \dim_RX$, and $\h(\fa, X)$.

\vspace{.1in}

\noindent{\bf Proposition 2.4.} For $X\in{\cal D}_+(R)$ the
following hold.

\begin{verse}

(a) $\cd(\fa, X)\le\cd(\fa, R)-\inf X$;

(b) $\cd(\fa, X)\le \dim_RX$.

\end{verse}
In addition, if $X\in{\cal D}_+^f(R)$ then

(c) $\h(\fa, X)\le\cd(\fa, X)$.

\vspace{.1in}

\noindent{\it Proof.} (a) and (b) are consequences of (1.0.1) and
the Grothendieck vanishing theorem on local cohomology. Part (c)
follows from the well--known fact that $\h(\fa, N)\le\cd(\fa, N)$
for any finite $R$--module $N$.\hfill$\square$

\vspace{.2in}

\noindent{\bf Proposition 2.5.} If $X\in{\cal D}_+(R)$, then

$$\cd(\fa, X)\ge\sup\{\cd(\fa, R/\fp)-\inf X_{\fp}|\fp\in\Spec
R\}.$$ Moreover, equality hold if $X\in{\cal D}_+^f(R)$.
\vspace{.1in}

\noindent{\it Proof.} We may assume that $\H(X)\neq 0$. For
$\fp\in\Supp X$, take $\ell=\inf X_{\fp}$ so that $\fp\in\Supp
\H_{\ell}(X)$. Thus we have $$\cd(\fa, R/\fp)-\inf
X_{\fp}\le\cd(\fa, \H_{\ell}(X))-\ell$$ which implies the
inequality.

Now assume $X\in{\cal D}_+^f(R)$. Let $\ell\in\Bbb Z$ such that
$\H_{\ell}(X)\neq 0$. By Theorem 1.3,
$\cd(\fa,\H_{\ell}(X))-\ell=\cd(\fa, R/\fq)-\ell$, for some
$\fq\in\Supp\H_{\ell}(X)$, which implies $\cd(\fa,
\H_{\ell}(X))-\ell\le\cd(\fa, R/\fq)-\inf X_{\fq}$, the result
follows. \hfill$\square$

\vspace{.2in}

\noindent{\bf Proposition 2.6.} If $X\in{\cal D}^f_+(R)$ and $M$
is a finite $R$--module, then
$$\cd(\fa, X\otimes^{\mathbf L}_RM)=\cd(\fa, X\otimes^{\mathbf
L}_RR/\Ann_R(M))\le\cd(\fa, X).$$

\vspace{.1in}

\noindent{\it Proof.} We have $X\otimes^{\mathbf L}_RM\in{\cal
D}^f_+(R)$ so, by Proposition 2.5,

\[ \begin{array}{rl}
\cd(\fa, X\otimes_R^{\mathbf L}M) &\, = \sup\{\cd(\fa, R/\fp)-\inf(X\otimes_R^{\mathbf L}M)_{\fp}|\fp\in\Spec R\}\\
&\, =\sup\{\cd(\fa, R/\fp)-\inf X_{\fp}|\fp\in\Supp X\cap\Supp M\}\\
&\, \le\sup\{\cd(\fa, R/\fp)-\inf X_{\fp}|\fp\in\Supp X\}\\
&\, =\cd(\fa, X).
\end{array}\]

To prove the first equality, we proceed as follows.
$$\cd(\fa, X\otimes_R^{\mathbf L}M)=\sup\{\cd(\fa, R/\fp)-\inf(X\otimes_R^{\mathbf L}M)_{\fp}|\fp\in\Spec
R\}.$$ By [{\bf Fo1}; Lemma 16.28], we have
$\inf(X_{\fp}\otimes^{\mathbf L}_{R_{\fp}}M_{\fp})=\inf
(X_{\fp}\otimes^{\mathbf
L}_{R_{\fp}}R_{\fp}/\Ann_{R_{\fp}}M_{\fp})$.\hfill$\square$

\vspace{.2in}

In [{\bf A}], Apassov defined the weak annihilator of a complex
$X\in{\cal D}(R)$ to be the intersection of the annihilators of
all the homology modules of $X$ and denoted by $\Ann_RX$.

The following result compares cohomological dimension of $X$ with
that of $R/\Ann_RX$. It is shown, in particular when $X\in{\cal
D}^f_+(R)$, that $\cd(\fa, X)$ lies in the interval $$[\cd(\fa,
R/\Ann_R X)-\sup X, ~\cd(\fa, R/\Ann_R X)-\inf X]$$ of length
$\amp X=\sup X-\inf X$.

\vspace{.1in}

\noindent{\bf Theorem 2.7.} Let $X\in{\cal D}_+(R)$. Then
$$\cd(\fa, R/\Ann_R X)-\sup X\le\cd(\fa, X)\le\cd(\fa, R/\Ann_R X)-\inf X.$$

\vspace{.1in}

\noindent{\it Proof.} By using Proposition 2.4, we can choose
$\ell\in\Bbb Z$ such that $\cd(\fa,
X)=\cd(\fa,\H_{\ell}(X))-\ell$. By Lemma 2.1, $\cd(\fa,
\H_{\ell}(X))\le\cd(\fa, R/\Ann_R X)$. As $\ell\ge\inf X$, the
right hand side inequality follows.

By Theorem 1.3, $\cd(\fa, R/\Ann_R X)=\cd(\fa, R/\fp)$ for some
prime ideal $\fp\supseteq\Ann_R X$. Hence
\[ \begin{array}{rl}
\cd(\fa, R/\Ann_R X) &\, \le \cd(\fa, R/\fp)-\inf X_{\fp}+\sup X\\
&\, \le\cd(\fa, X)+\sup X.
\end{array}\]
The last inequality follows from Proposition 2.5, which implies
the assertion. \hfill$\square$

\vspace{.2in}

The Theorem 2.9 expresses the cohomological dimension of
$X\otimes^{\mathbf L}_RY$ with the cohomological dimensions of the
tensor product of the homology modules of $X$ and $Y$. But first
we bring the following auxiliary result.

\vspace{.1in}

\noindent{\bf Proposition 2.8.} If $M$ and $N$ are finite
$R$--modules, then
$$\cd(\fa, M\otimes_R^{\mathbf L}N)=\cd(\fa, M\otimes_RN).$$

\vspace{.1in}

\noindent{\it Proof.} We have $M\otimes_R^{\mathbf L}N\in{\cal
D}^f_+(R)$, by [{\bf Fo1}; 7.28 and 7.31], and that, by
Proposition 2.5,

\[ \begin{array}{rl}
\cd(\fa, M\otimes_R^{\mathbf L}N) &\, = \sup\{\cd(\fa, R/\fp)-\inf(M\otimes_R^{\mathbf L}N)_{\fp}|\fp\in\Spec R\}\\
&\, =\sup\{\cd(\fa, R/\fp)-\inf(M_{\fp}\otimes_{R_{\fp}}^{\mathbf
L}N_{\fp})|\fp\in\Spec R\}\\
&\, =\sup\{\cd(\fa, R/\fp)-\inf M_{\fp}-\inf N_{\fp}|\fp\in\Spec
R\}
\end{array}\]

\noindent and the assertion follows by Theorem 1.3.\hfill$\square$

\vspace{.2in}

\noindent{\bf Theorem 2.9.} Let $X,~ Y\in{\cal D}^f_+(R)$. Then
$$\cd(\fa, X\otimes^{\mathbf L}_RY)=\sup\{\cd(\fa,
\H_{t}(X)\otimes_R\H_{\ell}(Y))-\ell-t|t, \ell\in\Bbb Z\}.$$

\vspace{.1in}

\noindent{\it Proof.} We have $X\otimes^{\mathbf L}_RY\in{\cal
D}^f_+(R)$, c.f. [{\bf Fo1}; 7.28 and 7.31]. By using [{\bf Fo1};
16.28] and Proposition 2.5, the following equalities hold.

\[ \begin{array}{rl}
\cd(\fa, X\otimes_R^{\mathbf L}Y) &\, = \sup\{\cd(\fa, R/\fp)-\inf(X\otimes_R^{\mathbf L}Y)_{\fp}|\fp\in\Spec R\}\\
&\, =\sup\{\cd(\fa, R/\fp)-\inf(X_{\fp}\otimes_{R_{\fp}}^{\mathbf
L}Y_{\fp})|\fp\in\Spec R\}\\
&\, =\sup\{\cd(\fa, R/\fp)-\inf_{\ell\in\Bbb Z}
(\inf(X_{\fp}\otimes_{R_{\fp}}^{\mathbf
L}\H_{\ell}(Y_{\fp}))+\ell|\fp\in\Spec R\}\\
&\, =\sup\{\cd(\fa, R/\fp)-\inf(X_{\fp}\otimes_{R_{\fp}}^{\mathbf
L}\H_{\ell}(Y_{\fp}))-\ell|\fp\in\Spec R, ~\ell\in\Bbb Z\}\\
&\, =\sup\{\sup\{\cd(\fa, R/\fp)-\inf(X\otimes_{R}^{\mathbf
L}\H_{\ell}(Y))_{\fp}|\fp\in\Spec R\}-\ell|\ell\in\Bbb Z\}\\
&\, =\sup\{\cd(\fa, X\otimes^{\mathbf
L}_R\H_{\ell}(Y))-\ell|\ell\in\Bbb Z\}.
\end{array}\]

\noindent Now the assertion holds by iterating the above technique
and using Proposition 2.7.\hfill$\square$

\vspace{.2in}

\noindent{\bf Corollary 2.10.} Let $X,~ Y\in{\cal D}^f_+(R)$. Then
$$\cd(\fa, X\otimes^{\mathbf L}_RY)\le\cd(\fa, X)-\inf Y.$$

\vspace{.1in}

\noindent{\it Proof.} It follows from Theorem 2.9 and Proposition
2.6.\hfill$\square$

\vspace{.3in}

\noindent{\bf 3. Cohomological dimension and the right derived
section functor}

\vspace{.2in}

For a complex $X\in{\cal D}_+(R)$, ${\mathbf R}\G_{\fa}(X)$, the
right derived section functor with support in $\V(\fa)$ applied to
the complex $X$, has been studied by several authors, e.g. [{\bf
AJL}], [{\bf FI}], [{\bf Fo2}], [{\bf Fr}], [{\bf L}], [{\bf
Sch}], and [{\bf Y}], whom have studied the invariant
$\sup{\mathbf R}\G_{\fa}(X)$. In this section we are interested in
$\inf{\mathbf R}\G_{\fa}(X)$.

The following result is a new form of [{\bf Fo2}; Proposition 2.5]
and the proof is similar, which H.-B. Foxby clarified to us in a
private discussion.

\vspace{.1in}

\noindent{\bf Lemma 3.1.} Let $X,~Y\in{\cal D}_+(R)$. Then

$$\inf(X\otimes^{\mathbf L}_RY)\ge\inf\{\inf(X\otimes^{\mathbf
L}_R\H_{\ell}(Y))+\ell|\ell\in\Bbb Z\}.$$ Moreover, the equality
holds if $X\in{\cal D}_+^f(R)$.\hfill$\square$

\vspace{.2in}

Now we are ready to prove the main result of this section.

\vspace{.1in}

\noindent{\bf Theorem 3.2.} If $X\in{\cal D}_+(R)$, then
$-\inf{\mathbf R}\G_{\fa}(X)\le\cd(\fa, X)$. Moreover, the
equality holds if $X\in{\cal D}^f_+(R)$.

\vspace{.1in}

\noindent{\it Proof.} By [{\bf Sch}; Proposition 3.2], there is a
functorial isomorphism $${\mathbf R}\G_{\fa}(X)\cong{\mathbf
R}\G_{\fa}(R)\otimes^{\mathbf L}_RX.$$ Hence we have the
following.

\[ \begin{array}{rl}
 -\inf{\mathbf R}\G_{\fa}(X) &\, =-\inf({\mathbf
 R}\G_{\fa}(R)\otimes^{\mathbf L}_RX)\\
&\, \le -\inf\{\inf({\mathbf
 R}\G_{\fa}(R)\otimes^{\mathbf L}_R\H_{\ell}(X))+\ell|\ell\in\Bbb
Z\}\\
&\, =-\inf\{-\cd(\fa, \H_{\ell}(X))+\ell|\ell\in\Bbb
Z\}\\
&\, =\sup\{\cd(\fa, \H_{\ell}(X))-\ell|\ell\in\Bbb Z\}\\
&\, =\cd(\fa, X).
\end{array}\]

\noindent The inequality follows by Lemma 3.1 and the second
equality by [{\bf Sch}; Proposition 3.2].

\noindent Now assume $X\in{\cal D}^f_+(R)$. Then the inequality
becomes equality by Lemma 3.1.\hfill$\square$

\vspace{.2in}

To present the final result we remind the notion of $\fa$--$\depth
X$, which is denoted by $\depth(\fa, X)$, and defined by
$$\depth(\fa, X)=-\sup{\mathbf R}\Hom_R(R/\fa, X).$$

In [{\bf FI}; Theorem 2.1], Foxby and Iyengar show that $\depth
(\fa, X)=-\sup{\mathbf R}\G_{\fa}(X)$.

\vspace{.1in}

\noindent{\bf Theorem 3.3.} Let $R$ be local, $X\in{\cal D}_+(R)$
and that $\depth(\fa, X)<\infty$. Then
$$\depth(\fa, X)\le\cd(\fa, X).$$

\vspace{.1in}

\noindent{\it Proof.} One has

\[ \begin{array}{rl}
 \depth(\fa, X) &\, =-\sup{\mathbf
 R}\G_{\fa}(X)\\
&\, \le -\inf{\mathbf
 R}\G_{\fa}(X)\\
&\, \le \cd(\fa, X).
\end{array}\]

Now the assertion holds.\hfill$\square$

\vspace{.2in}

We end this paper with two questions.

1. For any $X\in{\cal D}(R)$ and $\fp\in\Supp X$ we have the
inequality of dimensions $\dim X_\fp+\dim R/\fp\le \dim X$, cf.
[{\bf Fo1}; 16.16]. Is there any similar inequality for
cohomological dimensions?

2. If $X\in{\cal D}_+^f(R)$ is not homologically trivial and
$Y\in{\cal D}_b^f(R)$ then by the Intersection Theorem we have the
following inequality $\dim Y\le\dim (X\otimes^{\mathbf L}_RY)+\pd
X$. Is there any similar inequality for the cohomological
dimensions, such as $$\cd(\fa, Y)\le\cd(\fa,X\otimes^{\mathbf
L}_RY)+\pd X?$$

\vspace{.3in}

\noindent {\large\bf Acknowledgments.} The authors would like to
thank H.-B. Foxby, University of Copenhagen, for his invaluable
help, specially for posting the latest version of the reference
[{\bf Fo1}]. The research of the first author was partially
supported by a grant from IPM, and that of the second author was
partially supported by a grant from University of Tehran.
\vspace{.3in}

\baselineskip=16pt

\begin{center}
\large {\bf References}
\end{center}
\vspace{.2in}

\begin{verse}

[A] D. Apassov, {\em Annihilating complexes of modules}, Math.
Scand. {\bf 84} (1999), 11--22

[AJL] L. Alonso Tarr$\acute{\mathrm{i}}$o, A.
Jerem$\acute{\mathrm{i}}$as L$\acute{\mathrm{o}}$pez, J. Lipman,
{\em Local homology and cohomology on schemes}, Ann. Sci.
$\acute{\mathrm{E}}$cole Norm. Sup. (4) {\bf 30} (1997), 1--39.

[DNT]  K. Divaani-Aazar, R. Naghipour, M. Tousi, {\em
Cohomological dimension of certain algebraic varieties}, Proc.
Amer. Math. Soc. 130 (2002), no. 12, 3537--3544.

[F] G. Falting, {\em $\ddot{U}$ber lokale Kohomologiegruppen hoher
Ordnung}, J. Reine Angew. Math. {\bf 313} (1980), 43--51.

[G] A. Grothendieck, {\em Cohomologie locale des faisceaux
coh$\acute{e}$rents et th$\acute{e}$or$\grave{e}$mes de Lefschetz
locaux et globaux $(SGA$ $2)$}, Augment$\acute{\mathrm{e}}$
$\acute{\mathrm{d}}$un expos$\acute{\mathrm{e}}$ par
Mich$\grave{e}$le Raynaud. S$\acute{\mathrm{e}}$minaire de
G$\acute{\mathrm{e}}$om$\acute{\mathrm{e}}$trie
Alg$\acute{\mathrm{e}}$brique du Bois-Marie, 1962. Advanced
Studies in Pure Mathematics, Vol. 2. North-Holland Publishing Co.,
Amsterdam; Masson $\&$ Cie, $\acute{\mathrm{E}}$diteur, Paris,
1968.

[Ha1] R. Hartshorne, {\em  Residues and duality}, Lecture Notes in
Math., {\bf 20}, Springer Verlag, 1971.

[Ha2] R. Hartshorne, {\em Cohomological dimension of algebraic
varieties}, Ann. of Math. (2) {\bf 88} (1968), 403--450.

[Hu] C. Huneke, {\em Problems on local cohomology}, Free
resolutions in commutative algebra and algebraic geometry
(Sundance, UT, 1990), 93--108, Res. Notes Math., 2, Jones and
Bartlett, Boston, MA, 1992.

[HL] C. Huneke, G. Lyubeznik, {\em On the vanishing of local
cohomology modules} Invent. Math. {\bf 102} (1990), no. 1, 73--93.

[Fo1] H.-B. Foxby, {\em  Hyperhomological algebra and commutative
algebra}, Notes in preparation.

[Fo2] H.-B. Foxby, {\em Bounded complexes of flat modules}, J.
Pure Appl. Algebra {\bf 15} (1979), 149--172.

[FI] H.-B. Foxby, S. Iyengar, {\em Depth and amplitude for
unbounded complexes}, to appear in Contemporary Math.

[Fr] A. Frankild, {\em Vanishing of local homology}, to appear in
Math. Z.

[GM] J. P. C. Greenlees, J. P. May, {\em Derived functors of
$I$-adic completion and local homology}, J. Algebra {\bf 149}
(1992), 438--453.

[I] S. Iyengar, {\em Depth for complexes, and intersection
theorems}, Math. Z. {\bf 230} (1999), 545--567.

[L] J. Lipman, {\em Lectures on local cohomology and duality.
Local cohomology and its applications} (Guanajuato, 1999), 39--89,
Lecture Notes in Pure and Appl. Math., 226, Dekker, New York,
2002.

[Sch] P. Schenzel, {\em Proregular sequences, local cohomology,
and completions}, preprint.

[Y] S. Yassemi, {\em Generalized section functors}, J. Pure Appl.
Algebra {\bf 95} (1994), 103--119.

\end{verse}

E-mail addresses:

\hspace{1.2in} dibaeimt@ipm.ir

\hspace{1.2in} yassemi@ipm.ir

\end{document}